\documentclass[10.5pt]{article}
\usepackage{amsfonts,mathrsfs}

\def\sqr#1#2{{\vcenter{\vbox{\hrule height.#2pt
              \hbox{\vrule width.#2pt height#1pt \kern#1pt \vrule width.#2pt}
              \hrule height.#2pt}}}}
\def\signed #1{{\unskip\nobreak\hfil\penalty50
              \hskip2em\hbox{}\nobreak\hfil#1
              \parfillskip=0pt \finalhyphendemerits=0 \par}}
\def\endpf{\signed {$\sqr69$}}
\def\3n{\negthinspace \negthinspace \negthinspace }
\def\2n{\negthinspace \negthinspace }
\def\1n{\negthinspace }

\def\dbE{\mathbb{E}}
\def\dbF{\mathbb{F}}

\def\dbN{\mathbb{N}}

\def\dbP{\mathbb{P}}

\def\dbR{\mathbb{R}}

\def\dbY{\mathbb{Y}}
\def\dbZ{\mathbb{Z}}

\def\sG{\mathscr{G}}

\def\sK{\mathscr{K}}
\def\sL{\mathscr{L}}
\def\sM{\mathscr{M}}

\def\sP{\mathscr{P}}

\def\sX{\mathscr{X}}
\def\sY{\mathscr{Y}}
\def\sZ{\mathscr{Z}}

\def\={\buildrel \triangle \over =}

\def\ds{\displaystyle}

\def\ns{\noalign{\ss}}

\def\b{\beta}
\def\g{\gamma}
\def\d{\delta}
\def\e{\varepsilon}

\def\l{\lambda}

\def\si{\sigma}

\def\f{\varphi}

\def\Th{\Theta}
\def\L{\Lambda}

\def\O{\Omega}

\def\cF{{\cal F}}

\def\cK{{\cal K}}
\def\cL{{\cal L}}

\def\ss{\smallskip}
\def\ms{\medskip}

\def\q{\quad}
\def\qq{\qquad}
\def\hb{\hbox}

\def\esssup{\mathop{\rm esssup}}

\def\h{\widehat}
\def\wt{\widetilde}
\def\cd{\cdot}

\def\({\Big (}
\def\){\Big )}
\def\[{\Big[}
\def\]{\Big]}
\def\bde{\begin{definition}}
\def\ede{\end{definition}}
\def\be{\begin{equation}}
\def\bel{\begin{equation}\label}
\def\ee{\end{equation}}
\def\bt{\begin{theorem}}
\def\et{\end{theorem}}
\def\bc{\begin{corollary}}
\def\ec{\end{corollary}}
\def\bl{\begin{lemma}}
\def\el{\end{lemma}}
\def\bp{\begin{proposition}}
\def\ep{\end{proposition}}
\def\bas{\begin{assumption}}
\def\eas{\end{assumption}}
\def\br{\begin{remark}}
\def\er{\end{remark}}
\def\ba{\begin{array}}
\def\ea{\end{array}}
\def\ed{\end{document}}

\def\square#1{\vbox{\hrule\hbox{\vrule height#1%
     \kern#1\vrule}\hrule}}
\def\rectangle#1#2{\vbox{\hrule\hbox{\vrule height#1%
     \kern#2\vrule}\hrule}}

\font\tenbb=msbm10 \font\sevenbb=msbm7 \font\fivebb=msbm5

\newfam\bbfam
\scriptscriptfont\bbfam=\fivebb \textfont\bbfam=\tenbb
\scriptfont\bbfam=\sevenbb

\newtheorem{lemma}{Lemma}[section]
\newtheorem{remark}{Remark}[section]

\newtheorem{theorem}{Theorem}[section]
\newtheorem{corollary}{Corollary}[section]

\newtheorem{definition}{Definition}[section]
\newtheorem{proposition}{Proposition}[section]
\newtheorem{assumption}{Assumption}[section]

\makeatletter
   
   \@addtoreset{equation}{section}
\makeatother

\usepackage{amsfonts}
\usepackage{latexsym}
\usepackage{amsmath}
\usepackage{amssymb}
\usepackage{color}

\begin{document}

\title{Uniqueness of equilibrium strategies in dynamic mean-variance problems with random coefficients\footnote{The research was supported by the NSF of China under grant 11231007, 11401404 and 11471231.}}

\author{Tianxiao Wang \footnote{School of
Mathematics, Sichuan University, Chengdu, P. R. China. Email:wtxiao2014@scu.edu.cn.}}

\maketitle

\begin{abstract}
This paper is concerned with the uniqueness issue of open-loop equilibrium investment strategies of dynamic mean-variance portfolio selection problems with random coefficients. A unified method is developed to treat both the problems with deterministic risk-free return rate, state-dependent risk aversion, and that with full random coefficients, constant risk aversion. To do so, some new necessity conditions for the existence of equilibrium investment strategies are established, which considerably extends the analogue in \cite{Hu-Jin-Zhou-2017} with distinctive ideas.
Some new interesting facts are revealed. For example, if risk-free return rate is random, it is shown that there exists a unique open-loop equilibrium investment strategy relying on initial wealth, even when risk aversion is merely a constant but not state-dependent.

\end{abstract}

\ms


\ms


\section{Introduction}

The classical Markowitz's mean-variance portfolio selection problem is well-known in modern
investment portfolio theory. Due to practical requirement, it is natural to extend the original single period framework into multi-period scenario. If so, the
strategy at one moment may not be optimal at next moment. This phenomenon is the so-called time inconsistency.

In order to chase optimal choice, one has to keep changing the investment strategy in a naive way all the time.
On the other hand, the time consistency of policies is also fundamental in many situations. Therefore, the notion of time consistent equilibrium control/strategy has attracted much attention in the literature. Since mean-variance portfolio selection problems have some inherent connection with time inconsistent stochastic linear quadratic (SLQ) optimal control problems, we first take a revisit of the later.

As to time inconsistent stochastic linear quadratic problems, there are two types of equilibrium controls: open-loop equilibrium controls and closed-loop equilibrium controls. The first notion was introduced and investigated carefully in \cite{Hu-Jin-Zhou-2012}, \cite{Hu-Jin-Zhou-2017} with deterministic coefficients. The second one was considered in \cite{Bjork-Murgoci-2014} and \cite{Yong-2017} with two independent approaches. At this very moment, one may ask: what are the differences between these two notions? How are they relate to traditional optimal controls? If the time inconsistency in linear quadratic problems disappears, it was proved in \cite{Wang-2018-1} that open-loop equilibrium controls, which are characterized by first-order, second-order necessary optimality conditions, are strictly weaker than open-loop optimal controls, while closed-loop equilibrium controls reduce to closed-loop optimal controls. Related topics on time inconsistent SLQ problems can also be found in \cite{Djehiche-Huang-2016}, \cite{Wang-Wu-2016}, \cite{Wei-Yong-Yu-2017} and the references therein.

We return back to our mean-variance problems, where above two notions are referred as open-loop equilibrium (investment) strategies and closed-loop equilibrium (investment) strategies. The existence and uniqueness of the former notion was given in \cite{Hu-Jin-Zhou-2012}, \cite{Hu-Jin-Zhou-2017} when partial involved coefficients are random. Later, \cite{Wei-Wang-2017} extended it into the general asset-liability management problem with full random coefficients. As to closed-loop equilibrium investment strategy, it was firstly discussed in \cite{Basak-Chabakauri-2010} under the Markovian framework. Since the equilibrium strategy in \cite{Basak-Chabakauri-2010} does not depend on initial wealth, the authors in \cite{Bjork-Murgoci-Zhou-2012} introduced one class of mean-variance problems with state dependent risk aversion parameter, and discussed the existence of the associated equilibrium investment strategy. Later, \cite{Huang-Li-Wang-2017} developed a unified method to treat the problems in both \cite{Basak-Chabakauri-2010} and \cite{Bjork-Murgoci-Zhou-2012}, and obtained the uniqueness of closed-loop equilibrium investment strategies. Just recently, \cite{Wang-2018-2} extended the conclusions in \cite{Basak-Chabakauri-2010} into the general non-Markovian setting. One interesting fact is the obtained equilibrium strategies are allowed to rely on initial wealth, even though the risk aversion parameter is a constant.
We refer to \cite{Cui-Xu-Zeng-2016}, \cite{Czichowsky-2013}, \cite{We-et-al-2013}, \cite{Zeng-Li-2011}, etc., for more related investigations.

In this paper, we discuss the uniqueness of open-loop equilibrium investment strategy when all the coefficients are random. As we mentioned above, the uniqueness of open-loop equilibrium investment strategy was given in \cite{Hu-Jin-Zhou-2017}, where the risk-free return rate is deterministic. In \cite{Wei-Wang-2017}, the existence of open-loop equilibrium investment strategy was given with full random coefficients, while the uniqueness issue was absent. Therefore, the aim of our study is to fill this gap. On the other hand, as was shown in \cite{Wang-2018-1}, open-loop equilibrium controls (which is comparable with our open-loop equilibrium investment strategy) are normally weaker than open-loop optimal controls, even in classical time consistent setting. When discussing a weak notion, uniqueness is especially important from the theoretical point of view.


To study the uniqueness, we first obtain a necessary condition for open-loop equilibrium investment strategy with random coefficients. In contrast with \cite{Hu-Jin-Zhou-2017}, the randomness of risk-free return rate leads to the appearance of a unbounded stochastic process. As a result, the approach developed in \cite{Hu-Jin-Zhou-2012} fails to work here (see Remark \ref{Remark-Existing-fails-here}). To get round this difficult, we establish some proper convergence tricks (see Lemma \ref{Lemma-1-1}, Lemma \ref{One-lemma-e}) which have independent interest.
Using these necessity conclusions, we give a unified uniqueness treatment on mean-variance problems with deterministic risk-free return rate, state-dependent risk aversion, and that with full random coefficients, constant risk aversion. As to the former problem, our method is different from that in \cite{Hu-Jin-Zhou-2017}, and is simpler. As to the later one, there are some new interesting facts arising, see Subsection 4.1.

The rest of this paper is organized as follows. In Section 2, some notations, spaces are introduced and the financial problem is formulated. Section 3 is aim to provide new necessary conditions for the existence of open-loop equilibrium investment strategies. Section 4 is devoted to treating the uniqueness of equilibrium investment strategy. In Section 5, some concluding remarks are present.

%
%
%

\section{Preliminary notations and model formulation}

Through out this paper, $(\Omega, \cF,
\dbP, \{\cF_t\}_{t\geq0})$ is a filtered complete probability space, $\{W(t),
t\geq 0\}$ is a $\dbF:=\{\mathcal{F}_t\}_{t\geq 0}$-adapted one-dimensional Brownian
motion.

For
$H:=\dbR^n,\dbR^{n\times m}$, etc., $0\le s<t\le T$, $l\in\dbN$, we define
$$\ba{ll}
\ns\ds L^2_{\cF_t}(\Omega;H):=\Big\{X:\Omega\to H\bigm| X\hb{
is $\cF_t$ measurable,}\ \dbE|X|^2<\infty\Big\},\\
\ns\ds L^2_{\dbF}(s,t;H):=\Big\{X:[s,t]\times\Omega\to H\bigm|X(\cd)\hb{
is  measurable and}\\
\ns\ds\qq\qq\qq\q \hb{$\dbF$-adapted}, \ ~\dbE\Big[\int_s^t|X(r)|^2dr\Big]<\infty\Big\},\\
\ns\ds L^2_{\dbF}\big(\O;C([s,t];H)\big):=\Big\{X:[s,t]\times\O\to
H\bigm|X(\cd)\hb{ is measurable, \ $\dbF$-adapted,}\\
\ns\ds\qq\qq\qq\qq\qq \hb{has continuous paths,}\ \dbE\(\sup_{r\in[s,t]}|X(r)|^2\)<\infty\Big\},\\
\ns\ds L^\infty_{\dbF}(\Omega;C([s,t];\dbR)):= \Big\{X:[s,t]\times\O\to
H\bigm|X(\cd)\hb{ is measurable, $\dbF$-adapted, }\\
\ns\ds\qq\qq\qq\qq\qq \hb{ has continuous paths,}\ \esssup_{\omega\in\Omega}\sup_{r\in[s,t]}|X(r)|^2<\infty\Big\},\\
\ns\ds L^l_{\dbF}(\Omega;L^2(s,t;\dbR)):=\Big\{X:[s,t]\times\Omega\to H\bigm|X(\cd)\hb{
is  measurable and $\dbF$-adapted}\\
\ns\ds\qq\qq\qq \qq \q \ ~\dbE\Big[\int_s^t|X(r)|^2dr\Big]^{\frac l 2}<\infty\Big\}.
\ea$$

We consider a financial market where two assets are traded continuously on $[0,T]$. Suppose the price of bond evolves as
\begin{eqnarray*}
\left\{\begin{array}{rl}
dS_0(s) & \!\!\!= r(s)S_0(s)ds,\;\; s \in [0, T],\\
 S_0(0) & \!\!\!= s_0>0,
\end{array}\right.
\end{eqnarray*}
and the risky asset is described by
\begin{eqnarray*}
\left\{\begin{array}{rl} dS(s) & \!\!\!=
S(s)\Big\{b(s)ds+ \si (s)dW(s)\Big\},\;\;
s \in [0, T],\\
 S(0) & \!\!\!= s_1>0.\;\;
\end{array}\right.
\end{eqnarray*}
Here $r >0$ is the risk-free return rate, $b $ is the expected return rate of risky asset, $\si $ is the corresponding volatility rate.

\ms

(H0) Suppose $r,$ $b $, $\sigma$ are bounded and $\dbF$-adapted processes, and there exists constant $\d>0$ such that $|\sigma|^2 \geq \delta.$
%

\ms

Given initial capital $x>0$, $\beta:=b-r$, $\theta :=\beta \si^{-1}$, the investor's wealth $X$ satisfies
\begin{equation}\label{wealth-equation}
\!\!\!\!\!\!\left\{\begin{array}{rl}
\!\!\!dX(s)  & \!\!\!=\big[r(s) X(s)  +\beta(s)  u(s) \big] ds+u(s)\si(s)  dW(s),  \\
 \!\!\!X(0) & \!\!\!=x,
\end{array}\right.
\end{equation}
where $u $ is the capital invested in the risky asset.

At anytime $t\in[0,T)$, the objective of a mean-variance portfolio choice
model is to choose an investment strategy to minimize
\begin{eqnarray}\label{cost-functional}
J(u(\cd);t,X(t))=\hb{Var}_{t}\big[X(T)\big]-[\gamma_1+\gamma_2X(t)]\dbE_{t}\big[X(T)\big],
\end{eqnarray}
where $\dbE_{t}[\cdot]:=\dbE[\cdot|\mathcal{F}_t],$
$\gamma_1\geq0$, $\g_2\geq0$ are constants satisfying $\g_1\g_2=0$.

In this paper, we consider the open-loop equilibrium (investment) strategy of above problem. To this end, for $t\in[0,T)$, $v\in L^2_{\cF_t}(\Omega;\dbR)$, $\e>0$, we define $u^{v,\e}:=u^*+ vI_{[t,t+\e]}$.
\bde\label{Definition-1}
 Given $x\in\dbR$, $ u^*\in L^2_{\dbF}(0,T;\dbR)$
is called an {\it open-loop equilibrium investment strategy} if for any $t\in[0,T)$, $(X^*,u^*)$ satisfying (\ref{wealth-equation}),
\bel{optimal-open}\lim_{\overline{\e\to0}}
{J(u^{v,\e}(\cd);t,X^*(t))-J\big(u^*(\cd)\big|_{[t,T]};t,X^*(t)\big)
\over\e}\ge0.
\ee
\ede

In order to give a clearer picture of $u^*$, we introduce the following notion.

\bde\label{Definition-2}
A pair of $(\Th^*,\f^* )\in L^p_{\dbF}(\Omega;L^2(0,T;\dbR))\times L^2_{\dbF}(0,T;\dbR)$, $p>2$, is called an open-loop equilibrium operator if for any $x\in\dbR$, $u^*:=(\Th^*X^*+\f^*)\in L^2_{\dbF}(0,T;\dbR)$, and $u^*$ is an open-loop equilibrium investment strategy associated with $x$, where
\bel{Equilibrium-state-closed-loop}\left\{\!\! \ba{ll}
\ns\ds dX^*(s)=\big[[r(s)+\b(s)\Th^*(s)]X^*(s)+\b(s)\f^*(s)\big]ds\\
\ns\ds\qq\qq +\si(s)\big[\Th^*(s)X^*(s)+\f^*(s)\big]
dW(s),\ \ s\in[0,T],\\
\ns\ds X^*(0)=x.
\ea\right.
\ee
\ede

Here $(\Th^*,\f^*)$ does not depend on initial wealth $x$, while $u^*$ indeed does. In addition, the integrability of $(\Th^*,\f^*)$ ensures that (\ref{Equilibrium-state-closed-loop}) admits a unique continuous strong solution $X^*$.

\begin{remark}
Let us look at the condition $\g_1\g_2=0$, which includes if and only if (i): $\g_1=\g_2=0$; (ii): $\g_1=0$, $\g_2>0$; (iii) $\g_1>0$, $\g_2=0.$

For the second case, if risk-free return rate $r$ is deterministic, the equilibrium strategy $u^*:=\Th^* X^*$ indeed exists (\cite{Hu-Jin-Zhou-2012}). Plugging it into (\ref{wealth-equation}), the equilibrium wealth process $X^*(\cd)>0$. Hence
\bel{Standard-form-1-remark}\ba{ll}
\ns\ds
\frac{J(u^*(\cd);t,X^*(t))}{\g_2 X^*(t)}=\frac{1}{\g_2 X^*(t)}\hb{Var}_{t}\big[X^*(T)\big]- \dbE_{t}\big[X^*(T)\big],
\ea
\ee
where $\g(x):=\frac {1}{\g_2 x}$ represents the state-dependent risk aversion.

For the third case, the equilibrium strategy of form $u^*:=\Th^* X^*+\f^*$ also exists, even when $r$ is random (\cite{Wei-Wang-2017}). Hence
\bel{Standard-form-2-remark}\ba{ll}
\ns\ds
\frac{J(u^*(\cd);t,X^*(t))}{\g_1}=\frac{1}{\g_1}\hb{Var}_{t}\big[X^*(T)\big]- \dbE_{t}\big[X^*(T)\big],
\ea
\ee
where constant $\frac {1}{\g_1}$ represents the constant risk aversion.

Notice that the right hand of both (\ref{Standard-form-1-remark}) and (\ref{Standard-form-2-remark}) are the original and standard forms to describe the cost functional of mean-variance problems.

We look at the case when $\g_1\g_2=0$ does not hold, i.e., $\g_1>0,$ $\g_2>0$. According to Remark \ref{Remark-negative} next, it is impossible to always keep
$$X^*(t)\geq0,\ \ \g_1+\g_2 X^*(t)>0, \ \ t\in[0,T).$$
Consequently, there is no way to equivalently construct a cost functional of form
$$\h J(u(\cd);t,X(t)):=f(x)\hb{Var}_t\big[X^*(T)\big]-\dbE_t \big[X^*(T)\big], \ \ f'\leq 0,$$
by $J(u(\cd);t,X(t))$ such that $u^*$ is an equilibrium investment strategy in the sense of Definition \ref{Definition-1}. That is to say, the original interpretation of mean-variance problem is lost. Based on it, we will not discuss this type of problem.
\end{remark}

In the sequel, $K$ is a generic constant which varies in different context. Moreover, we may suppress the reference to the time variable whenever necessary.

\ms

\section{Necessary conditions of open-loop equilibrium investment strategy}

In this section, we derive some necessary conditions for the existence of open-loop equilibrium investment strategies. They help us discuss the uniqueness issue in the next section.

We first look at a motivational case of $\g_2=0$.
Given proper $(\Th,\f)$, $s\in[0,T],$ we consider
\bel{Equations-for-P-i}\left\{\ba{ll}
\ns\ds dP_1(s)=-\Big\{2r(s)P_1(s)+(P_1(s)\b(s)+\L_1(s)\si(s))\Th(s)\Big\}ds+\L_1(s)dW(s),\\
\ns\ds dP_2(s)=- r(s) P_2(s) ds+\L_2(s)dW(s),\\
\ns\ds dP_3(s)=- \big[ r(s) P_3(s)+(P_3(s)\b(s)+\L_3(s)\si(s))\Th(s)\big]ds+
\L_3(s)dW(s),\\
\ns\ds dP_4(s)=- \big[P_3(s)\b(s)+\L_3(s)\si(s)\big]\f(s) ds+\L_4(s)dW(s),\\
\ns\ds dP_5(s)=-\big[r (s)P_5(s) +\big[P_1(s)\b(s)+\L_1(s)\si(s)\big]\f(s)  \big]ds+
\L_5(s)dW(s),\\
\ns\ds P_1(T)=2,\ \ P_2(T)=-2,\ \ P_3(T)=1, \ \ P_4(T)=0,\ \ P_5(T)=-\g_1.
\ea\right.\ee
By Appendix B of \cite{Wei-Wang-2017}, there exist $(\Th,\f)$ and $(P_i,\L_i)$, $1\leq i\leq 5$,
satisfying both (\ref{Equations-for-P-i}) and
\bel{Condition-special-case}\left\{\ba{ll}
\ns\ds  \b (P_1 +P_2 P_3)+\si (\L_1
+\L_2 P_3)+\si^2 P_1 \Th=0,\ \ a.s. \ \ a.e. \\
\ns\ds  \b (P_2 P_4 +P_5 )+\si (\L_5 +\L_2 P_4 )+\si^2 P_1 \f=0.\ \ a.s. \ \ a.e.
\ea\right.
\ee
Moreover, under proper conditions (see Proposition 3.9 of \cite{Wei-Wang-2017}), for any $x\in\dbR$, $u:=\Th X+\f$ is an open-loop equilibrium investment strategy. In other words, $(\Th,\f)$ is an open-loop equilibrium operator in the sense of Definition \ref{Definition-2}.

At this moment, it is natural to ask: will all the open-loop equilibrium operator be obtained through (\ref{Condition-special-case})? Is (\ref{Condition-special-case}) a necessary condition for the existence of $(\Th,\f)$? We will give positive answers even when $\g_2\neq0$.

At first, let us make the following hypothesis.

\ms

(H1) There exists $(\Th,\f)\in L^p_{\dbF}(\Omega;L^2(0,T;\dbR))\times L^2_{\dbF}(0,T;\dbR)$, $p>2$ such that for any $x\in\dbR$, $u:=(\Th X+\f)\in L^2_{\dbF} (0,T;\dbR )$, where
\bel{State-equation-Th-1-Th-2}\left\{\2n\ba{ll}
\ns\ds
dX(s)=\big[(r(s)+\b(s)\Th(s)) X(s) +\b (s)\f(s) \big]ds \\
\ns\ds\qq \qq +\big[ \si (s) \Th (s) X(s)+\si (s) \f(s) \big]dW(s), \ \ s\in[0,T],\\
\ns\ds  X(0)=x.\ea\right.\ee

\ms

Fix $t\in[0,T)$, $\e>0$, $v\in L^2_{\cF_t}(\Omega;\dbR)$, let $u^{\e}:= u+vI_{[t,t+\e]}$ and $ X^\e $ be described as
\bel{Perturbed-1.1-ep}\left\{\2n\ba{ll}
\ns\ds
dX^\e(s)=\big[ r(s) X^\e(s)+\b(s)u^\e(s)\big]ds +  \si(s) u^\e (s)dW(s),  \ \ s\in[0,T], \\
\ns\ds X^\e(0)=x_0.
\ea\right.\ee
As a result, for $s\in[0,T],$ $X_0^\e(s):=X^\e(s)-X(s)$ satisfies
\bel{Equation-of-s-X-1-epislon}\left\{\2n\ba{ll}
\ns\ds
dX_0^\e(s) =\big[ r (s)  X^{\e}_0(s)   +\b (s) vI_{[t,t+\e]}(s)  \big]ds +  \si (s) vI_{[t,t+\e]} (s) dW(s), \\
\ns\ds X^\e_0(0)=0.\ea\right.\ee
The following estimate of $X_0^\e$ is easy to see,
$$\ba{ll}
\ns\ds \dbE_t\sup_{r\in[t,t+\e]}|X^\e_0(r)|^2\leq K\e,\ \  a.s., \ \ t\in[0,T).
\ea
$$
For $t\in[0,T], $ we introduce
\bel{MF-BSDEs-*-1}\left\{\ba{ll}
\ns\ds Y(s,t)=  2X(T)- 2\dbE_tX(T)-\g_1
+\int_s^{T}r(u) Y (u,t)du -\int_s^{T}Z (u,t)dW(u),\ \ s\in[t,T],\\
\ns\ds Y^\e_0(s,t)=  X^\e_0(T)- \dbE_tX^\e_0(T)
+\int_s^{T}r(u) Y^\e_0 (u,t)du -\int_s^{T}Z^\e_0(u,t)dW(u),\ \ s\in[t,T],\\
\ns\ds Y_0(s)=1+\int_s^Tr(u)Y_0(u)du-\int_s^TZ_0(u)dW(u),\ \ s\in[0,T].
\ea\right.\ee
Under (H1), system (\ref{MF-BSDEs-*-1}) is solvable with
$$
\ba{ll}
\ns\ds (Y(\cd,t),Z(\cd,t)), (Y^{v,\e}_0(\cd,t),Z^{v,\e}_0(\cd,t)) \in L^2_{\dbF}(\Omega;C([t,T];\dbR))\times L^2_{\dbF}(t,T;\dbR),\\
\ns\ds (Y_0(\cd),Z_0(\cd))\in L^\infty_{\dbF}(\Omega;C([0,T];\dbR))\times L^2_{\dbF}(t,T;\dbR),\ \ \ t\in[0,T].
\ea
$$
Moreover, it is easy to see that $(P_2,\L_2)=-2(Y_0,Z_0)$.

\bl\label{Cost-functional-1}
Suppose (H0), (H1) hold. Then
\bel{difference-cost}\ba{ll}
\ns\ds J(u^\e(\cd);t,X(t))-J(u(\cd);t,X(t)) \\
\ns\ds = \dbE_t\int_t^{t+\e}\b(s)\big[Y_0^{\e}(s,t)+Y(s,t)-\g_2X(t)Y_0(s)\big]ds \cd v\\
\ns\ds\qq+\dbE_t\int_t^{t+\e}\si (s)\big[Z_0^{\e}(s,t)+Z(s,t)-\g_2X(t)Z_0(s)\big]ds \cd v.
\ea\ee
\el

\it Proof. \rm From the definition of $u^\e$, it is a direct calculation that
$$\ba{ll}
\ns\ds J(u^\e(\cd);t,X(t))-J(u(\cd);t,X(t)) \\
\ns\ds = \dbE_t\Big\{\big[2X(T)- 2\dbE_tX(T)-(\g_1+\g_2 X(t))\big] X_0^\e(T)\Big\}\\
\ns\ds\qq +\dbE_t \Big\{( X_0^\e(T)- \dbE_tX_0^\e(T))X_0^\e(T)\Big\}.
\ea$$
By It\^{o}'s formula,
$$\left\{\ba{ll}
\ns\ds \dbE_t\Big\{\big[X_0^\e(T)- \dbE_tX_0^\e(T)\big] X_0^\e(T)\Big\}=\dbE_t \int_t^{t+\e}(\b(s)Y_0^\e(s,t)+\si(s)Z_0^\e(s,t))ds \cd v,\\
\ns\ds \dbE_t\Big\{\big[2X(T)- 2\dbE_tX(T)-\g_1 \big] X_0^\e(T)\Big\}=\dbE_t \int_t^{t+\e}(\b(s)Y(s,t)+\si(s)Z(s,t))ds \cd v,\\
\ns\ds \dbE_t \g_2 X(t) X_0^{\e}(T)=\g_2X(t)\dbE_t\int_t^{t+\e}(Y_0(s)\b(s)+Z_0(s)\si(s))ds \cd v.
\ea\right.
$$
Then above (\ref{difference-cost}) is easy to see.\endpf

\ms

We first treat the terms of $(Y,Z)$ in Lemma \ref{Cost-functional-1}. For $t\in[0,T]$, suppose
\bel{Representation-for-Y-*}\ba{ll}
\ns\ds  Y(s,t)=P_1(s)X(s)+ P_2(s)\dbE_t\big[P_3(s)X(s)+P_4(s)\big]+ P_{5}(s), \  \ s\in[t,T],
\ea
\ee
where $P_i$ satisfies
$$\left\{\ba{ll}
\ns\ds d P_i(s)=\Pi_i(s)ds+\L_i(s)dW(s), \ \ s\in[0,T],\ \  i=1,2,3,4,5,\\
\ns\ds P_1(T)=2,\ \ P_2(T)=-2,\ \ P_3(T)=1, \ \ P_4(T)=0,\ \ P_5(T)=-\g_1.
\ea\right.
$$
Here $\Pi_i$, $i=1,2,3,4,5$, is to be determined. By It\^{o}'s formula,
$$\left\{\ba{ll}
\ns\ds
d P_1X =\big[P_1 (r+\b\Th)+\Pi_1+\L_1 \si\Th
\big]X ds + (P_1[ \si\Th X+ \si\f]+\L_1X) dW(s)\\
\ns\ds\qq\qq\qq
+ (P_1\b+\L_1\si )\f ds,\\
\ns\ds d  P_2\dbE_t[ P_3X ] =\Big\{\Pi_2\dbE_t[
P_3X ]+P_2\Big\{\dbE_t\big[ (P_3r+P_3\b\Th + \Pi _3+\1n \L_3 \si\Th)X \big]
\\
\ns\ds\qq\qq\qq\qq +\dbE_t\big[(P_3\b+\L_3\si)\f \big] \Big\}\Big\}ds+ \L_2\dbE_t[P_3X]dW(s),\\
\ns\ds d P_2\dbE_t P_4 =[\Pi_2\dbE_t P_4+P_2 \dbE_t \Pi_4 ]ds+\L_2
\dbE_t P_4 dW(s).
\ea\right.$$
As a result, one has
$$\ba{ll}
\ns\ds
\Big\{\big[P_1 (r+\b\Th)+\Pi_1+\L_1 \si\Th
\big]X+ (P_1\b+\L_1\si )\f \\
\ns\ds\q + \Pi_2\dbE_t[
P_3X ]+P_2\Big\{\dbE_t\big[ (rP_3+\b P_3\Th + \Pi _3+\1n \L_3 \si \Th)X \big]
\\
\ns\ds\q+\dbE_t\big[(\b P_3+\L_3\si)\f \big] \Big\}+ \Pi_2\dbE_t P_4+P_2 \dbE_t \Pi_4 +\Pi_5\Big\}ds\\
\ns\ds \q+ \Big\{\L_2\dbE_t[P_3X] + P_1[ \si\Th X+ \si\f] +\L_2
\dbE_t P_4+\L_5+\L_1 X\Big\} dW(s)\\
\ns\ds =d\big[P_1X+P_2\dbE_t(P_3X+P_4)+P_5\big]=dY\\
\ns\ds=\big[-r P_1 X-rP_2\dbE_t[P_3X+P_4]-r P_5\big]ds+Z(s,t)dW(s).\ea$$
Comparing the diffusion term, one has
$$\ba{ll}
\ns\ds Z=\L_2\dbE_t[P_3X] + P_1[ \si\Th X+ \si\f] +\L_2
\dbE_t P_4+\L_5+\L_1 X.
\ea
$$
For the drift term,
$$\ba{ll}
\ns\ds P_1 (r+\b\Th)+\Pi_1+\L_1\si\Th+rP_1=0,\\
\ns\ds rP_2+\Pi_2=0,\ \ P_3r+P_3\b\Th + \Pi _3+\1n \L_3 \si\Th=0, \\
\ns\ds (P_3\b +\L_3\si)\f+\Pi_4=0,\ \ (P_1\b+\L_1\si)\f+rP_5+\Pi_5=0.
\ea
$$
Consequently, we introduce five BSDEs in above (\ref{Equations-for-P-i}) and make the following assumption.

\ms

(H2) Suppose there exists $(\Th,\f)$ such that for $P_i$ in (\ref{Equations-for-P-i}), and some $  p\in (2,\infty),$ $p'\in(1,2)$,
\bel{ }\ba{ll}
\ns\ds (P_i,\L_i)\in L^{\infty}_{\dbF}(\O;C([0,T];\dbR ))\times L^{p}_{\dbF}(\O;L^2(0,T;\dbR )),\ \ i=1,2,3,  \\
\ns\ds   (P_j,\L_j)\in L^{p'}_{\dbF}(\O;C([0,T];\dbR ))\times L^{p'}_{\dbF}(\O;L^2(0,T;\dbR )),\ \ j=4,5.
\ea\ee

\ms

For $t\in[0,T]$, and $s\in[t,T], $ we define
\bel{ }\left\{\ba{ll}
\ns\ds  M(s,t):=P_1(s)X (s)+P_2(s)\dbE_t[ P_3(s)X (s)+P_4(s)]+P_5(s), \\
\ns\ds N(s,t):=P_1(s)\si(s)[\Th(s) X(s)+\f(s)] +\L_1(s)X(s)\\
\ns\ds\qq\qq + \L_2(s)\dbE_t[P_3(s)X(s)+P_4(s)]+
\L_5(s).
\ea\right.\ee
Notice that $(M^{d},N^{d})$ are well defined, where
$$M^{d}(s):=M(s,s),\ \ N^{d}(s):=N(s,s),\ \ s\in[0,T].$$

\bl\label{Lemma-1-1}
Suppose (H0), (H1), (H2) hold. Then
\bel{Equality-1}\ba{ll}
\ns\ds \dbP\big\{\omega,\ Y(s,t)=M(s,t),\ \ \forall s\in[t,T]\big\}=1,\\
\ns\ds
 \dbP\big\{\omega,\  Z(s,t)=N(s,t)\big\}=1,\ \  s\in[t,T].\ \  a.e.
\ea\ee
Moreover, for almost $t\in[0,T]$, there exists sequence $\{\e_{n}^{(t)}\}_{n\geq1}$ such that $\lim\limits_{n\rightarrow\infty}\e_{n}^{(t)}=0$ and
\bel{Convergence-1}\ba{ll}
\ns\ds \lim_{n\rightarrow\infty}\Big[ \frac {1}{\e_{n}^{(t)}}\dbE_t\int_{t}^{t+\e_{n}^{(t)}}\big[\b(s) M(s,t)+\si(s)N(s,t)\big]ds\Big]\cd v \\
\ns\ds =\big[\b(t) M(t,t)+\si(t)N(t,t)\big]v.
\ea\ee

\el
To prove it, we need one more result.

\bl\label{One-lemma-e}
If $f(\cd)\in L^p_{\dbF}(0,T;\dbR^m)$, $m\in\dbN$, $p\in(1,2]$. Then for almost $t\in[0,T]$, there exists a sequence $\{\e^{(t)}_{n}\}$ depending on $t$ satisfying
$$\ba{ll}
\ns\ds \lim_{n\rightarrow\infty}\e^{(t)}_n=0,\ \ \ \lim_{n\rightarrow\infty}\frac {1}{\e^{(t)}_n}\int_{t}^{t+\e_n^{(t)}}\dbE_t|f_i(s)-f_i(t)|ds=0,\ \ i=1,2,\cdots,m.\ \ a.s. %
\ea
$$
\el
\it Proof. \rm
Given $f(\cd)$, we define $\wt f(\cd)$ as
$$\wt f(\cd):= f(\cd)I_{[0,T]}(\cd)+0I_{(T,2T]}(\cd).
$$
For $h>0$, there exists $g(\cd)\in C_{\dbF}([0,2T];L^p(\Omega;\dbR^m))$ such that
$$\big\|\wt f(\cd)-g(\cd)\big\|_{L^p_{\dbF}(0,2T;\dbR^m)}\leq h.$$
Hence
$$\int_0^{ T}\!\!\int_{t}^{t+\e}\dbE|\wt f(s)-g(s)|_{\dbR^m} ^p dsdt\!=\!\int_0^{\e}\!\!\int_{u}^{T+u}\dbE|\wt f(v)-g(v)|_{\dbR^m}^p dvdu\leq h\e.$$
%
%
As to $g(\cd)$, there exists $\d(h)>0$ such that
$$
 \dbE|g(t_1)-g(t_2)|^p_{\dbR^m} \leq h, \  \ |t_1-t_2|\leq\d,\ \ t_1,\ t_2\in[0,2T].
$$
As a result,
$$\int_0^{ T}\int_{t}^{t+\e}\dbE|g(s)-g(t)|_{\dbR^m}^pdsdt\leq h\e.$$
%
%
%
%
To sum up, for small $\e>0$,
$$\lim\limits_{\e\rightarrow0}\Big[\frac 1 \e \int_0^{ T}\int_{t}^{t+\e}\dbE|\wt f(s)
-\wt f(t)|_{\dbR^m}^pdsdt\Big]=0.$$
Hence there exists $\{\e_n\}_{n\geq1}$ such that for almost $t\in[0,T]$,
$$\lim\limits_{\e_n\rightarrow0} \frac {1} {\e_n} \big\|\wt f(s)-\wt f(t)\big\| ^p_{L^p_{\dbF}(t,t+\e;\dbR^m)}=0.$$
Notice that
$$ \dbE\Big[ \frac {1} {\e_n}  \int_{t}^{t+\e_n}\dbE_t|\wt f(s)-\wt f(t)|_{\dbR^m} ds \Big]^p\leq   \frac {1} {\e_n}  \int_{t}^{t+\e_n}\dbE |\wt f(s)-\wt f(t) |_{\dbR^m} ^pds.$$
Hence for almost $t\in[0,T],$ one selects a subsequence $\{\e_n^{(t)}\}$, such that
$$
\lim\limits_{\e_n^{(t)}\rightarrow0}\frac {1} {\e_n^{(t)}}  \int_{t}^{t+\e_n^{(t)}}
\dbE_t\big|\wt f(s)-\wt f(t)\big|_{\dbR^m} ds=0. \ \  a.s.$$
The definition of $\wt f(\cd)$ leads to the desired conclusion.\endpf

\ms

\it Proof of Lemma \ref{Lemma-1-1}. \rm
To prove (\ref{Equality-1}), we first claim that there exists $l\in(1,2)$ such that
\bel{Integrability-result-M-N}\ba{ll}
\ns\ds (M(\cd,t),N(\cd,t))\in L^{l}_{\dbF}(\Omega,C([t,T];\dbR))\times L^{l}_{\dbF}(\O;L^2(0,T;\dbR )),\ \ t\in[0,T].
\ea
\ee
Moreover, using It\^{o}'s formula, for any $t\in[0,T]$, we know that $(M(\cd,t),N(\cd,t))$ satisfies the first BSDE in (\ref{MF-BSDEs-*-1}). Therefore, (\ref{Equality-1}) is implied by the uniqueness of BSDEs.

Actually, under (H1), $X\in L^2_{\dbF}(\Omega,C([0,T];\dbR))$. The result of $M(\cd,t)$ in (\ref{Integrability-result-M-N}) is then followed from (H2).

Notice that $P_1(\cd)$ is bounded, from (H1), $P_1 D(\Th X+\f)\in L^2_{\dbF}(0,T;\dbR)$.
Thanks to the integrability in (H2), there exists $l\in(1,2)$,
$$
\left\{\ba{ll}
\ns\ds \dbE\Big[\int_0^T|\L_1X|^2ds\Big]^{\frac p 2}\leq \dbE\Big\{\sup_{t\in[0,T]}|X(t)|^p \Big[\int_0^T|\L_1 |^2ds\Big]^{\frac p 2}\Big\}<\infty,\\
\ns\ds \dbE\Big[\int_t^T\big|\L_2\dbE_t[P_3X ]\big|^2ds\Big]^{\frac l 2} +
\dbE\Big[\int_t^T\big|\L_2\dbE_t[P_4]\big|^2ds\Big]^{\frac l 2}
<\infty.
\ea\right.
$$
Then the result of $N(\cd,t)$ in (\ref{Integrability-result-M-N}) is easy to see.

Next we prove (\ref{Convergence-1}). By the definitions of $(M,N)$, for $u:=\Th X +\f$, $t\in[0,T)$,
$$
\ba{ll}
\ns\ds  \dbE_t \int_t^{t+\e}(\b(s)M(s,t)+\si(s)N(s,t))ds\\
\ns\ds =\dbE_t \int_t^{t+\e}\big[(\b(s)P_1(s)+\si(s)\L_1(s))X(s)+\b(s)P_5(s)+\si(s)\L_5(s)  +\si^2(s)P_1(s)u(s)\big]ds\\
\ns\ds\qq +\dbE_t \int_t^{t+\e}(\b(s)P_2(s)+\si(s)\L_2(s))\dbE_t(P_3(s)X(s)+P_4(s))ds.
\ea
$$
Under (H2), it follows from Lemma \ref{One-lemma-e} that there exists $\{\e_{n}^{(t)}\}$ such that
\bel{Set-of-convergence}\left\{
\ba{ll}
\ns\ds  \lim_{n\rightarrow\infty}\Big[\frac {1}{\e_{n}^{(t)}}\dbE_t \int_t^{t+\e_{n}^{(t)}}\big[\b P_2 +\si\L_2\big]ds\Big]= \b(t)P_2(t)+\si(t)\L_2(t),\ \ a.s.\\
\ns\ds  \lim_{n\rightarrow\infty}\Big[\frac {1}{\e_{n}^{(t)}}\dbE_t \int_t^{t+\e_{n}^{(t)}}(\b P_2+\si\L_2)(P_3X+P_4)ds\Big]\\
\ns\ds =(\b(t) P_2(t)+\si(t)\L_2(t))(P_3(t)X(t)+P_4(t)),\  \ a.s. \\
\ns\ds \lim_{n\rightarrow\infty}\Big[\frac {1}{\e_{n}^{(t)}}\dbE_t \int_t^{t+\e_n^{(t)}}\big[(
\b P_1+\si \L_1)X+\b P_5+\si \L_5 +\si^2 P_1 u\big]ds\Big]\\
\ns\ds= (\b(t)P_1(t)+\si(t)\L_1(t))X(t)+\b(t)P_5(t)+\si(t)\L_5(t) +\si^2(t)P_1(t) u(t).\ \ a.s.
\ea\right.
\ee
For any fixed $t\in[0,T]$,
$$
\ba{ll}
\ns\ds \dbE_t\sup_{s\in[t,t+\e]}\dbE_t |P_3(s)X(s)|^2\leq \dbE_t\sup_{s\in[t,t+\e]}|P_3(s)|^2\cd \dbE_t\sup_{s\in[t,t+\e]}|X(s)|^2<\infty,\ \ a.s.
\ea
$$
hence by conditional dominated convergence theorem,
$$\ba{ll}
\ns\ds \lim_{\e\rightarrow0}\dbE_t\sup_{s\in[t,t+\e]}|\dbE_t(P_3(s)X(s))-P_3(s)X(s)|^2=0.\ \ a.s.
\ea
$$
Similarly one has
$$\ba{ll}
\ns\ds \lim_{\e\rightarrow0}\dbE_t\sup_{s\in[t,t+\e]}|\dbE_tP_4(s)-P_4(s)|^2=0.\ \ a.s.
\ea
$$
Consequently,
\bel{}\ba{ll}
\ns\ds \lim_{\e\rightarrow0}\dbE_t\sup_{s\in[t,t+\e]}|\dbE_t(P_3(s)X(s)+P_4(s))-P_3(s)X(s)-P_4(s)|^2=0.\ \ a.s.
\ea
\ee
Using the first equality in (\ref{Set-of-convergence}), we then derive
$$\ba{ll}
\ns\ds \lim_{n\rightarrow\infty}\Big[\frac {1}{\e_{n}^{(t)}}\dbE_t \int_t^{t+\e_{n}^{(t)}}|\b P_2+\si\L_2|\Big|\dbE_t(P_3X+P_4)-P_3X-P_4\Big|ds\Big]=0. \ \ a.s.
\ea
$$
This, together with the second equality in (\ref{Set-of-convergence}), leads to
\bel{Second-last-1}\ba{ll}
\ns\ds \lim_{n\rightarrow\infty}\Big[\frac {1}{\e_{n}^{(t)}}\dbE_t \int_t^{t+\e_{n}^{(t)}}(\b P_2+\si\L_2)\dbE_t(P_3X+P_4)ds\Big]\\
\ns\ds =(\b(t)P_2(t)+\si(t)\L_2(t))(P_3(t)X(t)+P_4(t)).\ \
\ea
\ee
Eventually, our conclusion is followed by (\ref{Second-last-1}) and the third equality in (\ref{Set-of-convergence}). \endpf

\ms

For $(X,u)\in L^2_{\dbF}(\Omega;C([0,T];\dbR))\times L^2_{\dbF}(0,T;\dbR)$ satisfying (\ref{wealth-equation}), we define
\bel{Definition-of-s-Y-Z}\ba{ll}
\ns\ds  \sY :=\sP_1 X + \sP_2 \dbE_t\big[\sP_3 X +\sP_4 \big]+ \sP_{5},\\
\ns\ds \sZ:=\sL_2\dbE_t[\sP_3X] + \sP_1 \si u +\sL_2
\dbE_t \sP_4+\sL_5+\sL_1 X,
\ea
\ee
where for $s\in[0,T]$, $1\leq i\leq 5,$ $\sP_i(s)$ satisfies
\bel{Second-order-system}\left\{\ba{ll}
\ns\ds d\sP_1(s)=-2 r(s)\sP_1(s) ds+\sL_1(s)dW(s),\\
\ns\ds d\sP_2(s)=- r(s) \sP_2(s) ds+\sL_2(s)dW(s),\\
\ns\ds d\sP_3(s)=- r(s) \sP_3(s)    ds+
\sL_3(s)dW(s),\\
\ns\ds d\sP_4(s)=-  (\b(s)\sP_3(s)+\si(s)\sL_3(s))u(s) ds+\sL_4(s)dW(s),\\
\ns\ds d\sP_5(s)=-\big[r (s)\sP_5(s) +( \b(s)\sP_1(s)+\si(s)\sL_1(s))u(s) \big]ds+
\L_5(s)dW(s),\\
\ns\ds \sP_1(T)=2,\ \ \sP_2(T)=-2,\ \ \sP_3(T)=1, \ \ \sP_4(T)=0,\ \ \sP_5(T)=-\g_1.
\ea\right.\ee
By (H0), for $u\in L^2_{\dbF}(0,T;\dbR)$, and $p'\in(1,2)$, $p\in(2,\infty)$, (\ref{Second-order-system}) is solvable with
\bel{ }\ba{ll}
\ns\ds (\sP_i,\sL_i)\in L^{\infty}_{\dbF}(\O;C([0,T];\dbR ))\times L^{p}_{\dbF}(\O;L^2(0,T;\dbR )),\ \ i=1,2,3,\\
\ns\ds   (\sP_j,\sL_j)\in L^{p'}_{\dbF}(\O;C([0,T];\dbR ))\times L^{p'}_{\dbF}(\O;L^2(0,T;\dbR )),\ \ j=4,5.
\ea\ee
Moreover, recalling (\ref{Equations-for-P-i}), one has $(P_2,\L_2)\equiv (\sP_2,\sL_2)$.
\begin{corollary}\label{corollary-1}
Suppose (H0) holds. Then
\bel{Equality-2}\ba{ll}
\ns\ds \dbP\big\{\omega,\ Y(s,t)=\sY(s,t),\ \ \forall s\in[t,T]\big\}=1,\\
\ns\ds
 \dbP\big\{\omega,\  Z(s,t)=\sZ(s,t)\big\}=1,\ \  s\in[t,T].\ \  a.e.
\ea\ee
Moreover, for almost $t\in[0,T]$, there exists $\{\e_{n}^{(t)}\}_{n\geq1}$ such that $\lim\limits_{n\rightarrow\infty}\e_{n}^{(t)}=0$, and
\bel{Convergence-2}\ba{ll}
\ns\ds \lim_{n\rightarrow\infty} \frac {1}{\e_{n}^{(t)}}\dbE_t\int_{t}^{t+\e_{n}^{(t)}}\big[\b(s) \sY(s,t)+\si(s)\sZ(s,t)\big]ds\cd v \\
\ns\ds =\big[\b(t) \sY(t,t)+\si(t)\sZ(t,t)\big]v.
\ea\ee
\end{corollary}

\it Proof. \rm In Lemma \ref{Lemma-1-1}, let $\Th\equiv0$. Then
$$u=\f,\ \ (M,N)\equiv (\sY,\sZ),\ \ (P_i,\L_i)\equiv (\sP_i,\sL_i),\ \ i=1,2,3,4,5.$$
Our conclusions are implied by Lemma \ref{Lemma-1-1}.\endpf

\ms

Now let us turn to deal with the term with respect to $(Y^{v,\e}_0(\cd,\cd),Z^{v,\e}_0(\cd,\cd))$ in (\ref{difference-cost}).
\begin{lemma}\label{Lemma-3}
Suppose (H0) holds, $\sP_1$ is in (\ref{Second-order-system}). Then
\bel{3.15}\ba{ll}
\ns\ds
\lim_{\e\rightarrow0}\Big[\frac{1}{\e}\dbE_t\int_t^{t+\e}\big[
\b(s) Y^{v,\e}_0(s,t)+\si(s) Z^{v,\e}_0(s,t)\big]
ds \cd v\Big] =\frac 1 2 \si^2(t)\sP_1(t)v^2. \ \ a.s.
\ea\ee
\end{lemma}

\it Proof. \rm Firstly let us look
at the following equations on $[t+\e,T]$,
$$\left\{\2n\ba{ll}
\ns\ds
X^{v,\e}_0(s)=X^{v,\e}_0(t+\e)+\int_{t+\e}^s r(u)X^{v,\e}_0(u)du,\\
\ns\ds Y^{v,\e}_0(s,t)= X^{v,\e}_0(T)- \dbE_t[X^{v,\e}_0(T)]+\int_s^T r(u) Y^{v,\e}_0(u,t) du-\int_s^TZ^{v,\e}_0(u,t)dW(u).
\ea\right.$$
Using similar tricks as in Lemma \ref{Lemma-1-1}, one can prove that
$$\ba{ll}
\ns\ds \dbP\Big\{  Y^{v,\e}_0(s,t)=\frac 1 2\sP_1(s)X^{v,\e}_0(s)-\frac 1 4\sP_2(s)\dbE_t[\sP_2(s)X^{v,\e}_0(s)],\ \ \forall s\in[ t+\e,T]\Big\}=1,\\
\ns\ds \dbP\Big\{ Z^{v,\e}_0(s,t)=\frac 1 2\sL_1 (s) X^{v,\e}_0(s)
-\frac 1 2\sL_2(s)\dbE_t[\sP_2(s)X^{v,\e}_0(s)]\Big\}=1,\ \ s\in[t+\e,T].\ \ a.e.
\ea
$$
We continue to study
$(X^{v,\e}_0(\cd),Y^{v,\e}_0(\cd,t),Z^{v,\e}_0(\cd,t))$ on
$[t,t+\e]$ as follows,
$$\left\{\2n\ba{ll}
\ns\ds
X^{v,\e}_0(s)=\int_t^{s}\big[r(u)X^{v,\e}_0(u)+\b(u)v\big]du+\int_t^s \si(u)v dW(u),
\q  s\in[t,t+\e],\\
\ns\ds Y^{v,\e}_0(s,t)=Y^{v,\e}_0(t+\e,t)+\int_s^{t+\e}
 r(u)Y^{v,\e}_0(u,t) du -\int_s^{t+\e}Z^{v,\e}_0(u,t)dW(u),\q
 s\in[t,t+\e]. \ea\right.$$
Similar as Lemma \ref{Lemma-1-1}, by introducing
$$\left\{\2n\ba{ll}
\ns\ds d\wt \cK(s)= \frac 1 2(\sP_2(s)\b(s)+\sL_2(s)\si(s))vds+\wt\cL(s)dW(s), \ \ s\in[0,T],\\
\ns\ds  d \h \cK(s)=-\Big[r(s)\h\cK(s)+\frac 1 2(\sP_1(s)\b(s)+\sL_1(s)\si(s))v\Big]ds+\h\cL(s)dW(s),\\
\ns\ds  \wt\cK(t+\e)=\h\cK(t+\e)=0,\ea\right.
$$
we conclude that
$$\ba{ll}
\ns\ds
\dbP\Big\{\!Y^{v,\e}_{0}(s,t)\!=\!\frac 1 2 \sP_1(s)X^{v,\e}_0(s)\!-\!\frac 1 2 \sP_2(s)\dbE_t\big[\frac 1 2 \sP_2(s)X^{v,\e}_0(s)\!-\!\wt \cK(s)\big]\!+\!\h\cK(s),\ \forall s\in[t,t+\e]\!\Big\}\!=\!1,
\\
\ns\ds
\dbP\Big\{\!Z^{v,\e}_0(s,t)\!=\!\frac 1 2 \sL_1(s) X^{v,\e}_0(s)
\!+\!\frac 1 2 \sP_1(s)\si(s) v \!-\!\frac 1 2 \sL_2(s)\dbE_t\big[\frac 1 2 \sP_2(s)X^{v,\e}_0(s)\!-\!\wt\cK(s)\big]\!+\!\h\cL(s)\!\Big\}\!=\!1,
\ea
$$
with almost $s\in[t,t+\e].$ We observe that
\bel{Two-estimates-in-lemma}\left\{\!\!\ba{ll}
\ns\ds \dbE_t\sup_{s\in[t,t+\e]}|\wt \cK(s)|^2\!+\!\dbE_t\int_t^{t+\e}|\wt\cL(s)|^2ds \leq K\e v^2\dbE_t\int_t^{t+\e}\big[|\sP_2(s)|^2\!+\!|\sL_2(s)|^2\big]ds, \  \\
\ns\ds \dbE_t\sup_{s\in[t,t+\e]}|\h \cK(s)|^2\!+\!\dbE_t\int_t^{t+\e}|\h \cL(s)|^2ds \leq K\e v^2\dbE_t\int_t^{t+\e}\big[|\sP_1(s)|^2\!+\!|\sL_1(s)|^2\big]ds. \ \\
\ea\right.
\ee
As a result, for any $s\in[t,t+\e),$
\bel{B-Y-1-D-Z-1}\ba{ll}
\ns\ds \b(s) Y_0^{v,\e}(s,t)+\si(s) Z^{v,\e}_0(s,t)\\
\ns\ds
=\frac 1 2\big[\b(s) \sP_1(s)+\si(s)\sL_1 (s)
\big]X^{v,\e}_0(s)+\b(s)\h\cK(s)+\si(s)\h\cL(s)\\
\ns\ds\q -\frac 1 2\big[\b(s) \sP_2(s)+\si(s) \sL_2(s)\big]\dbE_t[\h \cK(s)X^{v,\e}_0(s)-\wt \cK(s)]+\frac 1 2\si^2(s) \sP_1(s) v.
\ea\ee
Recall that $\dbE_t\Big[\sup\limits_{_{t\in[t,t+\e)}}|X^{v,\e}_0(s)|^2\Big]
=O(\e)$, we then see that,
$$\left\{\2n\ba{ll}
\ns\ds\frac{1}{\e}\dbE_t\int_t^{t+\e}
\big[\b(s)\sP_1(s)+\si(s) \sL_1(s)
\big]X^{v,\e}_0(s) ds=o(1),\\
\ns\ds\frac{1}{\e}\dbE_t\int_t^{t+\e}
\big[\b(s) \sP_2(s)+\si (s)\sL_2 (s)\big]\dbE_t[\h \cK(s)X^{v,\e}_0(s)]
ds=o(1).\ea\right.$$
From (\ref{Two-estimates-in-lemma}), one has
$$\ba{ll}
\ns\ds\frac{1}{\e}\dbE_t\int_t^{t+\e} \big[\b(s)\h \cK(s)+\si(s)\h\cL(s)\big]ds=o(1),\\
\ns\ds \frac{1}{\e}\dbE_t\int_t^{t+\e}
\big[\b(s) \sP_2(s)+\si (s)\sL_2 (s)\big]\dbE_t[\wt\cK(s) ]
ds=o(1).
\ea$$
To sum up, our desired conclusion is followed by
$$\ba{ll}
\ns\ds \lim_{\e\rightarrow0}\Big[\frac 1 \e \dbE_t\int_t^{t+\e}\big[\b(s)Y_0^{v,\e}(s,t)+\si(s)Z_0^{v,\e}(s,t)\big]ds\Big]=\frac 1 2 \si^2(t)\sP_1(t)v.
\ea$$
\endpf

\ms

We present the first main result of this section.

\bt\label{Theorem-closed-open}
If $(\Th^*,\f^*)$ is an open-loop equilibrium operator in sense of Definition \ref{Definition-2} such that (H0), (H1), (H2) hold associated with $(P_i^*,\L_i^*)$. Then
\bel{Necessity-closed-open}\left\{\ba{ll}
\ns\ds \sG_1^*:= \b (P_1^* +P_2^* P_3^* -Y_0 \g_2)+\si (\L_1^*
+\L_2^* P_3^* -Z_0 \g_2)+\si^2 P_1^* \Th^*=0, \ \ a.s. \ a.e. \\
\ns\ds \sG_2^* :=\b (P_2^* P_4^* +P_5^* )+\si (\L_5^* +\L_2^* P_4^* )+\si^2 P_1^* \f^*=0. \ \ a.s. \ a.e.
\ea\right.
\ee
\et

\it Proof. \rm For $(Y_0,Z_0)$ in (\ref{MF-BSDEs-*-1}), we choose $\{\e_n^{(t)}\}_{n\geq 1}$ in Lemma \ref{Lemma-1-1} such that
\bel{One-convergence-result-Y-0}
\ba{ll}
\ns\ds \lim_{n\rightarrow\infty}\Big[\frac {1}{\e_n^{(t)}}\dbE_t\int_t^{t+\e_n^{(t)}}\big[\b(s)Y_0(s)+\si(s)Z_0(s)\big]ds \Big]  = \b(t)Y_0(t)+\si(t)Z_0(t).
\ea
\ee
From Lemma \ref{Cost-functional-1}, Lemma \ref{Lemma-1-1}, Lemma \ref{Lemma-3}, for almost $t\in[0,T]$,
\bel{ }\ba{ll}
\ns\ds 0\leq \big[ \sG_1^*(t)X^*(t)+\sG_2^*(t)\big] v+\frac 1 2 \si^2(t)\sP_1(t) v^2,
\ea
\ee
where $\sG_1^*$, $\sG_2^*$ are defined in (\ref{Necessity-closed-open}), $X^*$ is in (\ref{Equilibrium-state-closed-loop}). For any $x\in\dbR$, the arbitrariness of $v$ implies
$$\si^2(t)\sP_1(t)v^2\geq0,\ \   \sG_1^*(t)X^*(t)+\sG_2^*(t)=0.$$
By the integrability of $\Th^*$, the following equation admits a unique strong solution,
\bel{}
\left\{\2n\ba{ll}
\ns\ds
d\sX^*(s)= (r(s)+\b(s)\Th^*(s)) \sX^*(s) ds + \si(s) \Th ^*(s) \sX^*(s) dW(s), \ \ s\in[0,T],\\
\ns\ds  \sX^*(0)=1.
\ea\right.\ee
Moreover, $\big[\sX^*\big]^{-1}$ exists and is continuous. Suppose $X_0^*$ is the equilibrium wealth process corresponding to $X^*(0)=0$. We have
$$\ba{ll}
\ns\ds \dbP\big\{\omega\in\Omega;X^*(t)-X_0^*(t)=\sX^*(t),\ \ \forall t\in[0,T]\big\}=1.
\ea
$$
The arbitrariness of $x\in\dbR$ indicates that $\sG_1^* \sX^*=0$, $\sG_2^*=0$. In addition, the existence of $\big[\sX^*\big]^{-1}$ implies that $\sG_1^*=0$ as well.  \endpf

\br\label{Remark-negative}
We point out one interesting fact. Suppose the conditions in Theorem \ref{Theorem-closed-open} hold and $\b P_5^*+\si\L_5^*\neq0$. Notice that the later hypothesis is easy to be fulfill when $\g_1\neq0$. We claim that $\f^*\neq0$. Otherwise, according to (\ref{Equations-for-P-i}), $(P_4^*,\L_4^*)=(0,0)$, and the second equality in (\ref{Necessity-closed-open}) becomes $\b P_5^*+\si\L_5^*=0$, which is a contradiction. In other words, the equilibrium investment strategy is the form of $u^*:=\Th^* X^*+\f^*$. Plugging it into wealth equation (\ref{wealth-equation}), we see that $X^*$ may not always be non-negative.

\er
%

%
%
%

We give the second main result of this section.

\bt\label{Theorem-open}
Suppose (H0) holds, $(X',u')$ satisfies (\ref{wealth-equation}), $u'$ is an open-loop equilibrium investment strategy, $(\sY',\sZ')$ is defined in (\ref{Definition-of-s-Y-Z}) corresponding to $u'$. Then
\bel{Necessity-open}\ba{ll}
\ns\ds \b(t)\sY'(t,t)+\si(t)\sZ'(t,t)-\g_2X'(t)\big[ \b(t) Y_0(t)+\si(t)Z_0(t)\big]
=0.\ \ a.s. \  a.e.
\ea
\ee
\et

\it Proof. \rm  Given equilibrium strategy $u'$, from Lemma \ref{Cost-functional-1}, Corollary \ref{corollary-1}, Lemma \ref{Lemma-3}, and above (\ref{One-convergence-result-Y-0}), we have the following inequality for almost $t\in[0,T]$,
$$\ba{ll}
\ns\ds 0\leq \big[\b(t)(\sY'(t,t)+Y_0(t)\g_2X'(t))-\si(t)(\sZ'(t,t)+Z_0(t)\g_2X'(t))\big] v+\frac 1 2 \si^2(t)\sP_1(t) v^2,
\ea
$$
 By the arbitrariness of $v\in L^2_{\cF_t}(\Omega;\dbR)$, one has
$\si^2(t)\sP_1(t)v^2\geq0$ and above result. \endpf

\br\label{Remark-Existing-fails-here}
If coefficient $r$ is deterministic, an equivalent, yet less explicit form of (\ref{Necessity-open}) was obtained in Theorem 3.5 of \cite{Hu-Jin-Zhou-2017}. Nevertheless, when $r$ is random, $\L^*_2$ arise and is unbounded. Moreover, the tricks of Proposition 3.3 in \cite{Hu-Jin-Zhou-2017} fail. This explains the introducing of Lemma \ref{One-lemma-e}, which helps us bypass the encountered difficulty. As a trade off, we only derive the necessary condition, and leave the sufficiency verification for future study. More interestingly, it is just the necessity equilibrium condition that will help us solve the uniqueness issue next.
\er

\section{The uniqueness of open-loop equilibrium controls}

In this section, we discuss the uniqueness problem of open-loop equilibrium controls. To begin with, we make the assumption as follows.

\ms

(H3) For some $(\Th^*,\f^*)$, (\ref{Equations-for-P-i}) admits solution $(P_i^*,\L_i^*)$ such that $\si^2 P_1^*\geq\d>0$ and
\bel{One-important-fact}\ba{ll}
\ns\ds P_1^*+P_2^*P_3^*=0,\ \ \ \L_1^*+P_3^*\L_2^*+P_2^*\L_3^*=0,\ \ \int_0^{\cd}\L_1^*(s)dW(s) \text{ is BMO-martingale.}
\ea
\ee
\bt\label{Uniqueness-issue}
Suppose (H0) holds and $u^*:=\Th^*X^*+\f^*$ is an equilibrium control satisfying (H1), (H2), (H3). If $u'$ is another open-loop equilibrium control, then $u'=u^*$.
\et

\it Proof. \rm Given $(\sY',\sZ')$ in (\ref{Definition-of-s-Y-Z}) corresponding with equilibrium pair $(X',u')$, we define
$$
\left\{\ba{ll}
\ns\ds \bar Y':=\sY'-\big[P_1^*X'+P_2^*\dbE_t(P_3^*X'+P_4^*)+P_5^*\big]\equiv \sM_1+P_2^*\dbE_t \sM_2,\\
\ns\ds \bar Z':=\sZ'-\big\{P_1^*\si u' + \L_1^*X'+\L_2^*\dbE_t(P_3^*X'+P_4^*)+\L_5^*\big\} \equiv \L_2^*\dbE_t\sM_2+\sM_3,
\ea\right.
$$
where
\bel{}\ba{ll}
\ns\ds \sM_1:= (\sP_1-P_1^*)X'+(\sP_5-P_5^*), \ \ \sM_2:=(\sP_3-P_3^*)X'+\sP_4-P_4^*,\\
\ns\ds \sM_3:=(\sP_1-P_1^*)\si u'+(\sL_1-\L_1^*)X'+\sL_5-\L_5^*.
\ea\ee
Notice that in above we also use the fact that $(\sP_2,\sL_2)=(P_2^*,\L_2^*)$.
By It\^{o}'s formula,
$$
\left\{\ba{ll}
\ns\ds d \big[P_1^* X'\big]=\Big\{-r P_1^*X'  -(P_1^*\b+\L_1^*\si)\big[\Th^* X'-u'\big] \Big\}ds+\big[\L_1^*X'+P_1^* \si u' \big]dW(s),\\
\ns\ds d \big[P_3^* X' \big]=  -\big[P_3^*\b+\L_3^*\si\big]\big[\Th^*X'-u'\big]ds+\big[\L_3^*X'+P_3^* \si u'\big]dW(s),\\
\ns\ds d \big[P_2^*\dbE_t(P_3^*X')\big]\!=\!- \!\big[r P_2^* \dbE_t(P_3^*X')\!+\!
P_2^*\dbE_t\big[(P_3^*\b+\L_3^*\si)(\Th^*X'-u')\big] \big]ds\!+\!\L_2^*\dbE_t(P_3^*X')dW(s),\\
\ns\ds d \big[P_2^*\dbE_tP_4^*\big] =-\Big[P_2^*\dbE_t\big[(P_3^*\b+\L_3^*\si)\f\big] +r P_2^* \dbE_tP_4^*\Big]ds+\L_2^*\dbE_tP_4^*dW(s),
\ea\right.
$$
and
$$
\left\{\ba{ll}
\ns\ds d \big[\sP_1 X'\big]=\Big\{-r\sP_1 X' +(\sP_1\b+\sL_1\si) u' \Big\}ds+\big[\sL_1X'+\sP_1 \si u' \big]dW(s),\\
\ns\ds d \big[\sP_3 X' \big]=  \big[\sP_3\b+\sL_3\si\big] u' ds+\big[\sL_3 X'+\sP_3\si u'\big]dW(s),\\
\ns\ds d \big[\sP_2\dbE_t(\sP_3X')\big]=\Big\{- r \sP_2 \dbE_t(\sP_3X')+
\sP_2\dbE_t\big[(\sP_3\b+\sL_3\si)u'\big] \Big\}ds+\sL_2\dbE_t(\sP_3X')dW(s),\\
\ns\ds d \big[\sP_2\dbE_t\sP_4\big] =-\Big[\sP_2\dbE_t\big[(\sP_3\b+\sL_3\si)u'\big] +r \sP_2 \dbE_t\sP_4\Big]ds+\sL_2 \dbE_t\sP_4 dW(s).
\ea\right.
$$
Then it is a direct calculation that
\bel{Equations-for-sM-1-2}\left\{\ba{ll}
\ns\ds d\sM_1=\Big[-r\sM_1+( P_1^*\b+\L_1^*\si)[\Th^*  X' +\f^* -u' ]\Big]ds+\sM_3dW(s),\\
\ns\ds d\sM_2=( P_3^* \b +\L_3^*\si)[\Th^*  X' +\f^* -u' ] ds+\sM_4 dW(s),\\
\ns\ds \sM_1(T)=\sM_2(T)=0,
\ea\right.\ee
where
$$\sM_4:=\big[(\sL_3-\L_3^*) X'+(\sP_3-P_3^*)\si u'+\sL_4-\L_4^*\big].$$
On the other hand, notice that
\bel{Expression-of-db-Y-Z}\left\{
\ba{ll}
\ns\ds
\dbY'(s):=\bar Y'(s,s)=\sM_1(s)+P_2^*(s)\sM_2(s), \\
\ns\ds \dbZ'(s):=\bar Z'(s,s)=\L_2^*(s)\sM_2(s)+\sM_3(s),\ \ s\in[0,T],
\ea\right.
\ee
are well defined.
From Theorem \ref{Theorem-open}, Theorem \ref{Theorem-closed-open}, a necessary condition for the equilibrium control $u'(\cd)$ is
\bel{Necessary-condition-application}
\ba{ll}
\ns\ds 0=\b(s)(\sY'(s,s)-\g_2 X'(s)Y_0(s))+\si(s)(\sZ'(s,s)-\g_2 X'(s)Z_0(s))\\
\ns\ds\ \ =\b(s) \dbY'(s)+\si(s)\dbZ'(s)+\Big[\b(s)(P_1^*(s)+P_2^*(s)P_3^*(s)-\g_2 Y_0(s))+\si(s)(\L_1^*(s)\\
\ns\ds\qq +\L_2^*(s)P_3^*(s)-\g_2Z_0(s))\Big]X'(s) +\Big[\b(s) (P_5^*(s) +P_2^*(s) P_4^*(s) )\\
\ns\ds\qq +\si(s) [ \L_5^*(s) +\L_2^*(s)  P_4^*(s) ]\Big]+\si^{2}(s) P_1^*(s) u' (s) \\
\ns\ds\ \ =[\b(s)\dbY'(s)+\si(s)\dbZ'(s)] -\si^2(s) P_1^*(s)(\Th^*(s) X'(s)+\f^*(s)-u'(s)). \ \ a.s.
\ea
\ee
Plugging it into (\ref{Equations-for-sM-1-2}), for $s\in[0,T],$ we arrive at
$$
\left\{\ba{ll}
\ns\ds d\sM_1 =\Big[-r\sM_1 +\frac{[P_1^*\b+\L_1^*\si]}{\si^2P_1^*}\big[\b\dbY'+\si\sZ'-\si P_2^*\sM_4\big]\Big]ds+\sM_3 dW(s),\  \\
\ns\ds d\sM_2 = \frac{(P_3^*\b+\L_3^*\si)}{\si^2P_1^*}
\big[\b \dbY' +\si\sZ'-\si P_2^*\sM_4\big] ds+\sM_4 dW(s),\\
\ns\ds \sM_1(T)=\sM_2(T)=0,
\ea\right.
$$
where
\bel{Definition-of-sZ'}\ba{ll}
\ns\ds \sZ':=\sM_3+\L_2^*\sM_2+P_2^*\sM_4.
\ea\ee
Using It\^{o} formula to $P_2^*\sM_2$, and recalling (\ref{One-important-fact}),
$$\ba{ll}
\ns\ds d \big[P_2^*\sM_2\big] =\Big[P_2^*(P_3^*\b+\L_3^*\si) (\si^2P_1^*)^{-1}\big[\b\dbY' +\si\sZ'-\si P_2^*\sM_4\big]\\
\ns\ds \qq\qq\q-rP_2^*\sM_2+\L_2^*\sM_4\Big]ds+[\L_2^*\sM_2+P_2^*\sM_4]dW(s)\\
\ns\ds\qq\qq=\Big[-\big(P_1^*\b+(\L_1^*+P_3^*\L_2^*)\si\big)(\si^2P_1^*)^{-1}\big[\b
\dbY' +\si\sZ'-\si P_2^*\sM_4\big]\\
\ns\ds \qq\qq\q-rP_2^*\sM_2+\L_2^*\sM_4\Big]ds+[\L_2^*\sM_2+P_2^*\sM_4]dW(s).
\ea
$$
Again thanks to (\ref{One-important-fact}), we obtain the following equation for $(\dbY',\sZ')$
$$\ba{ll}
\ns\ds d\dbY'=\Big[-r\dbY'+\L_2^*\sM_4-P_3^*\L_2^* \si^{-1}\big[P_1^*\big]^{-1}
\big[\b\dbY' +\si\sZ'-\si P_2^*\sM_4\big]\Big]ds+\sZ' dW(s)\\
\ns\ds\qq=\Big[-r\dbY'+P_2^*\L_2^* \si^{-1}\big[\b\dbY' +\si\sZ'\big]\Big]ds+\sZ' dW(s).
\ea
$$
Recall $\dbY'(T)=0$ and the BMO property of $\L_2^*$, we obtain that
$$\ba{ll}
\ns\ds \dbP\big\{\omega\in\Omega;\ \dbY'(s)=0,\ \forall s\in[0,T]\big\}=1,\ \  \dbP\big\{\omega\in\Omega;\ \sZ'(s)=0\big\}=1,\ \ s\in[0,T],
\ea
$$
with the help of Theorem 10 in \cite{Briand-Confortola-2008}. Notice that
$$\ba{ll}
\ns\ds -\frac{(P_3^*\b+\L_3^*\si)}{\si^2P_1^*} \si P_2^*  =\frac{\big[P_1^*\b +(\L_1^*+P_3^*\L_2^*)\si \big]}{\si P_1^*}  =\big[\frac{\b}{\si }+\frac{\L_1^*}{P_1^*}-\frac{\L_2^*}{P_2^*}\big].
\ea
$$
Hence
$$\left\{\ba{ll}
\ns\ds d\sM_2=\big[\frac{\b}{\si }+\frac{\L_1^*}{P_1^*}-\frac{\L_2^*}{P_2^*}\big]\sM_4 ds+\sM_4 dW(s), \ \ s\in[0,T],\\
\ns\ds \sM_2(T)=0.
\ea\right.
$$
Thanks to the BMO property of $\L_1^*$ in (\ref{One-important-fact}), and  Theorem 10 in \cite{Briand-Confortola-2008}, we have $(\sM_2,\sM_4)=(0,0)$. Putting it back to (\ref{Expression-of-db-Y-Z}), (\ref{Definition-of-sZ'}), one has $(\sM_1,\sM_3)=(0,0)$. Hence $\dbZ'=0.$ Recalling $\dbY'=0,$ and above (\ref{Necessary-condition-application}), one has $u'=\Th^* X'+\f^*.$ Our conclusion follows naturally. \endpf

\br\label{Remark-comparisons-HHZ-2017}
Similar as \cite{Hu-Jin-Zhou-2017}, we adopt the following tricks: introducing $(\bar Y',\bar Z')$, $(\dbY',\dbZ')$, using the previous necessary conditions, and proving $(\dbY',\dbZ')=(0,0)$. Nevertheless, when coefficient $r $ is random instead of deterministic, $\L_2^*\neq0$, $\bar Z'$ depends on parameter $t$, and the tricks in Theorem 5.2 of \cite{Hu-Jin-Zhou-2017} do not work any more.

To figure out this issue, we point out two important observations. In the first place, thanks to the definitions of $(\sY',\sZ')$ in (\ref{Definition-of-s-Y-Z}), we can further represent explicitly $(\bar Y',\bar Z')$, $(\dbY',\dbZ')$ by means of introduced $\sM_i$, $i=1,2,3$. In the second place, the equations satisfied by $\sM_i$ can fortunately be derived by direct calculations. By introducing $\sZ'$ in (\ref{Definition-of-sZ'}) and some subtle deductions, we end up with a one-dimensional linear BSDE of $(\dbY',\sZ')$, and get the desired conclusion.
\er

\subsection{Mean-variance problems with constant risk aversion}

In this part, we discuss the mean-variance problem with constant risk aversion $\g_1$, and $\g_2=0$.

We first introduce the corresponding Riccati system. To this end, for $s\in[0,T]$, we consider
\bel{Equations-for-P-i-special-2}\left\{\ba{ll}
\ns\ds dP_1^*=-\Big\{2rP_1^*+(P_1^*\b+\L_1^*\si)\Th^*\Big\}ds+\L_1^*dW(s),\\
\ns\ds dP_2^*=- rP_2^* ds+\L_2^*dW(s),\\
\ns\ds P_1^*(T)=2,\ \ P_2^*(T)=-2,\ \ \Th^*:=-\frac{\L_1^*}{\si P_1^*}+\frac{\L_2^*}{\si P_2^*}.
\ea\right.\ee
  If we define
$$P_3^*:=-\frac{P_1^*}{P_2^*},\ \ \L_3^*:= -\frac{\L_1^*}{P_2^*}+\frac{P_1^*\L_2^*}{[P_2^*]^2},$$
then
\bel{P_3^*-constant}\left\{\ba{ll}
\ns\ds dP_3^* =- \big[ r P_3^* +(P_3^* \b +\L_3^* \si )\Th^* \big]ds+
\L_3^* dW(s),\\
\ns\ds P_3^*(T)=1.
\ea\right.\ee
Moreover, $\Th^*=-\frac{\L_3^*}{P_3^*}$.
Using $(P_3^*,\L_3^*)$, $(P_1^*,\L_1^*)$, for $s\in[0,T]$, we look at
\bel{}\left\{\ba{ll}
\ns\ds %
  dP_4^*=- \big[P_3^*\b+\L_3^*\si\big]\f^* ds+\L_4^*dW(s),\\
\ns\ds dP_5^* =-\big[r P_5^* +\big[P_1^* \b +\L_1^*\si \big]\f^*  \big]ds+
\L_5^*dW(s),\\
\ns\ds  P_4^*(T)=0,\ \ P_5^*(T)=-\g_1,\\
\ns\ds \f^*:=-\frac{ \b (P_2^* P_4^* +P_5^* )+\si (\L_5^* +\L_2^* P_4^* )}{\si^2 P_1^*}.
\ea\right.\ee
We verify the assumption in Theorem \ref{Uniqueness-issue} one by one.

By Appendix B, C in \cite{Wei-Wang-2017}, there exists $(P_i^*,\L_i^*)$, $i=1,2,3,4,5$ satisfying (H2), and (H1) holds with previous defined $(\Th^*,\f^*)$.

According to the definitions of $(P_3^*,\L_3^*)$, the first two equality in (\ref{One-important-fact}) hold. In addition, from Remark 3.6 in \cite{Wei-Wang-2017}, we see that $\int_0^{\cd}\L_1^*(s)dW(s)$ is a BMO-martingale.

To sum up, by applying Theorem \ref{Uniqueness-issue}, we obtain the uniqueness of open-loop equilibrium investment strategy $u^*:=\Th^* X^*+\f^*$.

\br
We point out some interesting facts.

(1) If $r$ is deterministic, by the uniqueness of BSDE (\ref{P_3^*-constant}), see Appendix B of \cite{Wei-Wang-2017}, one has $\L_3^*=0$ which leads to $\Th^*=0.$ In other words, the appearance of $\Th^*$ in $u^*$ is determined by the randomness of $r$.

(2) Even though the risk aversion parameter is a constant, the equilibrium investment strategy $u^*$ is allowed to depend on initial wealth as along as $r$ keeps some randomness.

(3) If $\b$, $\si$ are deterministic, or even $\b=0$, $\Th^*$ still does not degenerate as long as $r$ is random.
However, if $r$ becomes deterministic, $\Th^*=0$, even though $\b$, $\si$ are random. These two aspects show the different roles of $r,$ $\b$, $\si$ in effecting $\Th^*$.

(4) By introducing
$$\left\{\ba{ll}
\ns\ds dM_1=-\Big[\sK_1 M_1+\sK_2 N_1\Big]ds+N_1dW(s),\\
\ns\ds dM_2=-\Big[\sK_3 M_1+\sK_4N_1+\sK_5 N_2\Big]ds+N_2dW(s),\\
\ns\ds M_1(T)=-\g_1,\ \ M_2(T)=0.
\ea\right.
$$
where $\sK_i$ are defined as,
$$\left\{\ba{ll}
\ns\ds \sK_1:=(r-\frac{\b \L^*_2}{\si P_2^*}), \ \ \sK_2:=-\frac{\L_2^*}{P_2^*},\ \  \sK_3:= \frac{[P_3^*\b+\L_3^*\si ]}{\si^2 P_2^* P_3^*}\b,\\
\ns\ds  \sK_4:=\frac{[P_3^*\b+\L_3^*\si ]}{\si^2 P_2^* P_3^*}\si, \ \ \sK_5:=-\frac{[P_3^*\b+\L_3^*\si]}{\si P_3^*},
\ea\right.
$$
we can rewrite $\f^*$ as
\bel{}\ba{ll}
\ns\ds \f^*=\frac{\b M_1+\si N_1}{\si^2 P_2^* P_3^*}-\frac{N_2}{\si P_3^*}
\ea
\ee
If $r$ is deterministic, it is a direct calculation that
\bel{}\ba{ll}
\ns\ds \Th^*=0,\ \ \f^*=\big[\frac{\b \g_1}{2\si^2}-\frac{N_2}{\si P_3^*}\big]e^{-\int_{\cd}^T r(s)ds}.
\ea
\ee

(5) Our uniqueness ensures that there is no more pair of $(\Th^*,\f^*)$ satisfying above properties.

\er
%

%
%
%
%

\subsection{Mean-variance problems with state-dependent risk aversion}

In this part, we consider the mean-variance problems when $\g_1=0$, and $r$ be deterministic.

At first, we consider the following system of equations,
\bel{Equations-for-P-i-special-one}\left\{\ba{ll}
\ns\ds dP_1^*(s)=-\Big\{2r(s)P_1^*(s)+(P_1^*(s)\b(s)+\L_1(s)\si(s))\Th^*(s)\Big\}ds+\L_1^*(s)dW(s),\\
\ns\ds dP_2^*(s)=- r(s) P_2^*(s) ds,\ \  s\in[0,T],\\
\ns\ds P_1^*(T)=2,\ \ P_2^*(T)=-2,
\ea\right.\ee
where
$$\Th^*:=-\frac{\frac 1 2 \b P_2^* \g_2-\si\L_1^*}{\si^2P_1^*}.$$
The solvability of $(P_1^*,\L_1^*)$ was discussed in Section 5 of \cite{Hu-Jin-Zhou-2012}.
If we define
$$P_3^*:=-\frac{P_1^*}{P_2^*},\ \ \L_3^*:=-\frac{\L_1^*}{P_2^*},$$
we have
\bel{}\left\{\ba{ll}
\ns\ds dP_3^*(s)=- \big[ r(s) P_3^*(s)+(P_3^*(s)\b(s)+\L_3^*(s)\si(s))\Th^*(s)\big]ds+
\L_3^*(s)dW(s),\\
\ns\ds P_3^*(T)=1.
\ea\right.\ee
According to \cite{Hu-Jin-Zhou-2012}, \cite{Hu-Jin-Zhou-2017}, assumptions (H1), (H2), (H3) hold with $u^*:=\Th^* X^*$, $\f^*=0$. From Theorem \ref{Uniqueness-issue}, the uniqueness is easy to see.

\br
Notice that our procedures in proving above uniqueness are distinctive and simpler than the analogue in Section 5 of \cite{Hu-Jin-Zhou-2017}. More details along this can be found in Remark \ref{Remark-comparisons-HHZ-2017}.
\er
\br
In our situation, if $\g_2=0,$ it follows from Subsection 4.1 that $\Th^*$ reduces to zero. Similar conclusion also happens when $\b=0$. These facts indicate the important role of $\g_2$, $\b$ in keeping the non-degeneration of $\Th^*$ or $u^*$.
\er
\br
Let us compare the deterministic $r$ and $\g_2$. If $\g_2\neq0$, $\Th^*$ may not equal to zero even when $r=0$. However, if $\g_2=0$, $\Th^*$ must equal to zero even when $r\neq0$. This shows that $\g_2$ is more essential than deterministic $r$.
\er
%

%
%
%
%

\ms

\section{Concluding remark}

In this paper, the uniqueness of open-loop equilibrium investment strategies for dynamic mean-variance portfolio selection problems is discussed. A unified method is proposed to treat mean-variance problems in two different settings and some interesting aspects are revealed as well. We emphasize that there is no essential difficulty to extend the conclusions here into the case when both investment strategy and Brownian motions are multi-dimensional.

As to the problems with random risk-free return rate and state-dependent risk aversion, the existence and uniqueness of open-loop equilibrium investment strategies are still under consideration. We hope to discuss it in future publications.

\end{document}